\documentclass[11pt]{amsart}
\usepackage{anyfontsize}
\usepackage{tikz-cd}
\usepackage[headings]{fullpage}
\usepackage{amsmath}
\usepackage{amsfonts}
\usepackage{latexsym}
\usepackage{amssymb}
\usepackage{placeins}

\newtheorem{proposition}{Proposition}[section]
\newtheorem{theorem}{Theorem}[section]

\newtheorem{lemma}{Lemma}[section]

\newtheorem{corollary}{Corollary}[section]

\def\homo{\mathop{\sf Hom}}
\def\ker{\mathop{\sf Ker}}

\def\dbydt{\sf \frac{d}{dt}}

\def\ext{\sf Ext}

\def\Cinf{\mathcal{C}^\infty}
\def\Dr{\mathcal{D}'}

\def\Sr{\mathcal{S}'}

\def\D{\mathcal{D}}

\def\C{\mathbb{C}}

\def\N{\mathbb{N}}
\def\R{\mathbb{R}}

\newcommand{\del}{\partial}

\newcommand{\pde}{[\partial_1, \ldots ,\partial_n]}
\newcommand{\diff}{[\sigma_1, \ldots ,\sigma_n]}

\newcommand{\F}{\mathcal{F}}

\begin{document}

\title[Canonical Controller]{The Canonical Controller for Distributed Systems}

\maketitle

\begin{center}
\author{Shiva Shankar \footnote{Department of Electrical Engineering, IIT Bombay, \tt{shunyashankar@gmail.com}}}
\end{center}

\begin{abstract} This paper generalises results of Willems-Trentelman, and van der Schaft, on achievable behaviours, to the case of  linear distributed systems defined by partial differential or difference equations. It shows that the `minimal' controller which achieves a particular subsystem is the canonical controller of van der Schaft, thereby answering the `open problem' of \cite{sc} in the setting of infinite dimensional and $n-D$ systems. This result is used to describe the collection of all linear subsystems of the electro-magnetic field, containing the vacuum solutions, that can be attained by suitable choices of electric charge and current density.
\end{abstract}

{\tiny  \hspace{7mm} Keywords: Achievable behaviours; Partial differential and difference equations; Maxwell equations.}

\hspace{.7cm} {\tiny Mathematics Subject Classification: 93B10,  93C20, 93B25} 

\section{introduction}

This paper generalises the results in \cite{wt} on achievable subsets of behaviours described by ordinary differential equations, to the setting of distributed systems described by partial differential or difference equations. Furthermore, it shows that there is a unique minimal controller which accomplishes a given restriction, whose structure is identical to that of the canonical controller of \cite{sc}. 

Recollect that by definition, a behaviour is the collection of all the evolutions, or trajectories, of a dynamical system. Potentially, any possible evolution could perhaps occur,  but the laws that the system must obey, restrict the actual occurance to a subset. These laws, if they are local, are described by differential or difference equations. For example, the components of the electric and magnetic fields in space and time, could a priori have been  arbitrary functions in $\Cinf(\R^4)$, but in fact must satisfy the Maxwell equations. In this interpretation, a law serves to restrict the possible evolutions of a system, and the collection of all the laws that it satisfies, defines the system. If a system's trajectories must be further restricted, according to some criteria (such as stability, or rapid decay at infinity), then one must impose further laws, in the form of a controller. This world view, initiated by Willems \cite{w}, does not require notions of inputs or outputs in its formulation, and is a vast generalization of the classical state space theory.

The problem addressed by Willems-Trentelman in \cite{wt} is to characterize all the `achievable behaviours' of  a dynamical system. In this  problem, the variables which describe the attributes of a system are of two kinds, those which  need to be controlled, denoted $w$, and the variables by means of which control is accomplished, denoted $c$. In the input-output paradigm of state space theory, or in the transfer function approach, $c$ would be the inputs to the system, and $w$ its outputs. The trajectories of $w$ and $c$ that can possibly occur, and the relationships between them, are described by the laws of the system. The problem now is to characterize those  subsets of trajectories to which it is possible to restrict the evolution of the $w$ variables, by imposing restrictions on the control variables. These are the achievable behaviours of the system. A complete characterization of these behaviours is obtained in \cite{wt} for linear time invariant systems described by ordinary differential equations. The work in \cite{sc} provides a partial generalization in a very general context. Furthermore, in it, van der Schaft establishes the existence of a `canonical' controller that accomplishes the task of restriction, 
and in the process  uncovers an instance of the Internal Model Principle in its description. 

This paper generalises all these results to distributed systems. The description of achievable subsystems is a faithful generalisation of the results in \cite{wt}. There are several conrollers which restrict a distributed system to a given subsystem, and amongst them, there is a minimal controller, which turns out to be the faithful generalisation of the canonical controller of \cite{sc}. 
 This paper is thus an answer to the `open problem' in \cite{sc}, in the setting of infinite dimensional and $n-D$ systems.

The paper concludes with an application of these results to the control of the electro-magnetic field in space and time. For instance, if the electric and magnetic fields are the variables that must be controlled by suitable choices of electric current and charge density, then the results here provide a complete characterisation of 
the possible subsystems of the electro-magnetic field that can be so achieved.
 
\section{achievable subsystems of a  distributed systems} 
 
Let $A$ be either the ring $\C[\partial] = \C\pde$ of partial differential operators on $\R^n$, or the ring $\C[\sigma] = \C\diff$ of partial difference operators on the subset $\N^n \subset \mathbb{Z}^n$ of lattice points with positive integral coordinates. The attributes of the systems that we study take values in the space $\Dr$ of distributions on $\R^n$, in the first case, and in the space $\C^{\N^n}$, of all complex valued functions on $\N^n$, in the second. An element of $\C[\partial]$ acts on an element in $\Dr$ by differentiation, and gives $\Dr$ the structure of a $\C[\partial]$-module. Similarly, the action of $\sigma_i$ on an $f \in \C^{\N^n}$ by shift in the $i$-th coordinate, namely $\sigma_i(f)(x_1,\ldots, x_i, \ldots, x_n) = f(x_1, \ldots, x_i+1, \ldots x_n)$, makes $\C^{\N^n}$ a $\C[\sigma]$-module. More generally, the attributes of the system lie in an $A$-submodule of $\Dr$ or $ \C^{\N^n}$, as the case may be, for instance the{ space $\Cinf$ of smooth functions in $\Dr$, or the submodule of bounded functions in $\C^{\N^n}$. We call this $A$-submodule `the space of signals', and denote it by $\F$.\\

Let $P$ be an $A$-submodule of $A^k$, $k \geqslant 1$, the free $A$-module of rank $k$. It is finitely generated, say by $p_1, \ldots, p_\ell$, where $p_i = (p_{i1}, \ldots, p_{ik})$, ~$1 \leqslant i \leqslant \ell$. This choice of generators defines the matrix operator
\[
\begin{array}{cccc}
P(\cdot): & \F^k  & \longrightarrow & \F^\ell\\
& f=(f_1, \ldots, f_k) & \mapsto &  (p_1f, \ldots, p_\ell f) ,
 \end{array}
\]
where $P(\cdot)$ is either the partial differential operator $P(\del)$, or the partial difference operator $P(\sigma)$, depending upon the choice of $A$, and  where $p_if = \sum_{j=1}^k p_{ij}f_j$, for all $i$. The distributed system defined by $P(\cdot)$, in the signal space $\F$, is the kernel $\ker_{\F}(P(\cdot))$ of the above operator.
In Willems' interpretation, the rows of $P(\cdot)$ are the laws that determine the system, and to say $P(\cdot)f = 0$ is to say that $f$ satisfies these laws.

This kernel, however, depends not on the choice of generators for $P$ which make up the rows of $P(\cdot)$, but only on the  submodule $P$. Indeed, by Malgrange \cite{m}, the above kernel is isomorphic to the $A$-module $\homo_{A}(A^k/P, ~\F)$ of all $A$-linear maps from the quotient module $A^k/P$ to the signal space $\F$. This isomorphism is given by the map
\begin{equation}
\begin{array}{ccc}
\ker_\F(P(\cdot))  & \longrightarrow & \homo_A(A^k/P, ~\F)  \\ 
 f=(f_1, \ldots, f_k) & \mapsto & \phi_f   ~,
 \end{array}
\end{equation}
where $\phi_f([e_i]) = f_i$, $1 \leqslant i \leqslant k$, and where $[e_1], \ldots ,[e_k]$ denote the images of the standard basis $e_1, \ldots ,e_k$ of $A^k$ in $A^k/P$. The inverse of this map is
\[
\begin{array}{ccc}
\homo_A(A^k/P,~ \F)  & \longrightarrow & \ker_\F (P(\cdot))\\
 \phi & \mapsto & (\phi([e_1]), \ldots , \phi([e_k]))  ~.
 \end{array}
\]
  
Hence, we denote this kernel  by $\ker_\F(P)$, and call it `the system defined by the kernel of $P$ in $\F$'. An element $f \in \ker_\F(P)$ is a trajectory of the system.  

Clearly, $P \subset P'$ implies that $\ker_\F(P') \subset \ker_\F(P)$. \\  

As explained earlier, the attributes $f$ of the system $\ker_\F(P) \subset \F^k$, are of two types, the control variables $f_c \in \F^{k_c}$, and the variables $f_w \in \F^{k_w}$ that are to be controlled, where $k_w + k_c = k$. We then write an $f \in \F^k$ as $(f_w, f_c) \in \F^{k_w+k_c}$. Correspondingly, we denote an element $p \in A^k$ by $(p_w, p_c) \in A^{k_w+k_c}$. There are several injections and surjections defined by this separation of variables, and we denote them as in the following split exact sequences:
\[
0 \rightarrow A^{k_w} \begin{array}{c} \stackrel{\iota_w}{\longrightarrow}\\\stackrel{\pi_w}{\longleftarrow} \end{array} A^{k_w+k_c} \begin{array}{c}\stackrel{\pi_c}{\longrightarrow}\\\stackrel{\iota_c}{\longleftarrow} \end{array} A^{k_c} \rightarrow 0 ~,
\]
and 
\[
0 \leftarrow \F^{k_w} \begin{array}{c} \stackrel{\pi_w}{\longleftarrow}\\ \stackrel{\iota_w}{\longrightarrow} \end{array} \F^{k_w+k_c} \begin{array}{c}\stackrel{\iota_c}{\longleftarrow}\\ \stackrel{\pi_c}{\longrightarrow} \end{array} \F^{k_c} \leftarrow 0 ~,
\]
where the second sequence is obtained from the first by applying the functor $\homo_A( - ,~\F)$ to it.

An $A$-submodule $P$ of $A^{k_w+k_c}$ then defines the $A$-submodules $\pi_w(P)$ and $\iota^{-1}_w(P)$ of $A^{k_w}$, and the submodules $\pi_c(P)$ and $\iota^{-1}_c(P)$ of $A^{k_c}$. Similarly, the system $\ker_\F(P) \subset \F^{k_w+k_c}$ defined by $P$, defines the $A$-submodules $\pi_w(\ker_\F(P))$ and $\iota^{-1}_w(\ker_\F(P))$ of $\F^{k_w}$, and $\pi_c(\ker_\F(P))$ and $\iota^{-1}_c(\ker_\F(P))$ of $\F^{k_c}$. \\

\noindent Remark 2.1: This notation is slightly different from the notation in \cite{sc}, for instance the $A$-module $\iota^{-1}_w(\ker_\F(P))$ above, is denoted by $\mathcal{P}_0$ there. The notation here is meant to emphasise the interchangeable roles of $w$ and $c$. This symmetry is again observed below, in the comment after Lemma 2.2.\\

In this terminology, we can state the control problem of this paper: \\
(i) Consider the set of trajectories in $\pi_w(\ker_\F(P))$. Suppose that we wish to restrict it to a subset consisting of only those trajectories which satisfy some criterion defined by the problem, and that this is to be achieved by restricting the control trajectories in $\pi_c(\ker_\F(P))$ to a subset.

\vspace{1.5mm}
What are the subsets of $\pi_w(\ker_\F(P))$ that can be achieved by this process? \\

In the behavioural paradigm of Willems explained above, a restriction of the control trajectories is  achieved by imposing additional laws that the control variables must satisfy. These additional laws constitute the controller. We can then ask:
\vspace{1.5mm}

(ii) What is the structure of the controller that accomplishes the control task of (i)? \\

The rest of this section is dedicated to answering these questions. 
\vspace{1.5mm}

We first establish relationships between the various $A$-submodules of $\F^{k_w}$ and $\F^{k_c}$ derived from $\ker_\F(P)$ above.

\begin{proposition} Let $\F$ be any signal space, and $P$ be an $A$-submodule of  $A^{k_w+k_c}$. Then the submodules of $\F^{k_w}$ described above satisfy
\[ 
\iota_w^{-1}({\ker}_\F(P)) = {\ker}_\F(\pi_w(P)) \subset \pi_w({\ker}_\F(P)) \subset {\ker}_\F(\iota_w^{-1}(P)) ~, \] 
and similarly for the corresponding submodules of the control signal space $\F^{k_c}$.
\end{proposition}
\noindent Proof: Set $k = k_w + k_c$. Applying the left exact functor $\homo_A( -, ~\F)$ to the exact sequence $A^k/P \stackrel{\pi_w}{\longrightarrow} A^{k_w}/\pi_w(P) \rightarrow 0$, implies that   $0 \rightarrow \homo_A(A^{k_w}/\pi_w(P), ~\F) \stackrel{\iota_w}{\longrightarrow} \homo_A(A^k/P, ~\F)$ is exact. The isomorphism of (1) proves the equality in the statement of the proposition. 

Similarly, the short exact sequence 
\[ 0 \rightarrow A^{k_w}/\iota_w^{-1}(P) \stackrel {\iota_w}{\longrightarrow} A^k/P  \stackrel{\pi_c}{\longrightarrow} A^{k_c}/\pi_c(P) \rightarrow 0
\]
yields the exact sequence 
\[ 0 \rightarrow {\ker}_\F(\pi_c(P)) \stackrel{\iota_c}{\longrightarrow} {\ker}_\F(P) \stackrel{\pi_w}{\longrightarrow} {\ker}_\F(\iota_w^{-1}(P)) ~,
\]
and hence the second inclusion of the statement. 

Finally, if $f_w$ is in $\iota_w^{-1}({\ker}_\F(P))$, then $(f_w,0)$ is in $\ker_\F(P)$, and hence $f_w$ is also in $\pi_w(\ker_\F(P))$.
\hspace*{\fill}$\square$\\

As we work in the category of systems that arise as kernels of differential or difference operators, the first problem we encounter is that the projection of a system need not always be a system. \\

\noindent Example 2.1: Let $A=\mathbb{C}[\dbydt]$, and let $\F = \D$, the space of compactly supported smooth functions on $\R$. Let $P \subset A^2$ be the cyclic submodule generated by $(-1, ~\dbydt)$, and let $\pi_w: A^2 \rightarrow A$ be the projection to the first factor.

Then, $\ker_\D(P) = \{({\dbydt}f, f)~|~ f \in \D \}$,  and $\pi_w(\ker_\D(P))  = \{{\dbydt}f ~|~ f \in \D \}$. If this image were the  kernel of a differential operator, say the kernel of  $p(\dbydt): \D \rightarrow \D$, then it would follow that the composition $p(\dbydt) \circ \dbydt  = 0$. As $A$ is a domain, this implies that $p(\dbydt) = 0$, and hence that ${\dbydt}: \D \rightarrow \D$ is surjective. This is a contradiction, as the image of $\dbydt$  consists of only those elements in $\D$ that integrate to 0 on $\R$.\\

We overcome this problem by restricting the choice of the signal space $\F$ to an injective $A$-module. Recollect that to say $\F$ is injective, is to say that $\homo_A( -, ~\F)$ is an exact functor. The celebrated Fundamental Principle of Palamadov and Malgrange asserts that $\Dr$, $\Cinf$, and  the space $\Sr$ of temperate distributions, are injective $\C[\partial]$-modules. Moreover, it is an elementary fact that the space $\C^{\N^n}$ is an injective $\C[\sigma]$-module.

Also recollect that an injective $A$-module $M$ is a cogenerator if $\homo_A(P, M)$ is nonzero whenever $P$ is nonzero. The $\C[\partial]$-modules $\Dr$ and $\Cinf$ are cogenerators, whereas $\Sr$ is not a cogenerator, for instance \cite{sst}. Again it is elementary that $\C^{\N^n}$ is a cogenerator as a $\C[\sigma]$-module.

It follows that if $\F$ is injective, and a cogenerator, then there is an inclusion reversing bijection between $A$-submodules $P$ of $A^k$ and systems $\ker_\F(P)$ in $\F^k$, \cite{sst}. 

\begin{proposition} Let the signal space $\F$ be an injective $A$-module. Then the projection of a system is also a system. Furthermore, 
\[\pi_w({\ker}_\F(P)) = {\ker}_\F(\iota_w^{-1}(P)) ,\]
i.e. the second inclusion of Proposition 2.1 is an equality. 

Similarly, $\pi_c({\ker}_\F(P)) = {\ker}_\F(\iota_c^{-1}(P))$.
\end{proposition} 
\noindent Proof: It suffices to observe that the short exact sequence in the proof of Proposition 2.1, now yields a short exact sequence upon applying the exact functor $\homo_A( - ,~\F)$, and hence that $\pi_w: \ker_\F(P) \rightarrow \ker_\F(\iota_w^{-1}(P))$ is surjective. \hspace*{\fill}$\square$\\
 
\noindent Remark 2.2: When the signal space $\F$ is not an injective $A$-module, for example the space $\D$ of Example 2.1, then the obstruction to the above projection being a kernel lies in ${\ext}_A^1(A^{k_c}/\pi_c(P), ~\F)$ (see for instance \cite{sg}).\\

Hence, we assume for the rest of the paper that the space of signals $\F$ is an injective $A$-module.

We now state again the control problem that we study for such signal spaces: \\

{\em (i) Given a distributed system $B = \ker_\F(P) \subset \F^{k_w+k_c}$, it defines two other systems by projection, $B_w = \pi_w(B) = \ker_\F(\iota_w^{-1}(P)) \subset \F^{k_w}$, the system that is to be controlled, and $B_c = \pi_c(B) = \ker_\F(\iota_c^{-1}(P)) \subset \F^{k_c}$, the controller. The problem is to restrict $B_w$ to a desired subsystem by restricting the behaviour of the controller $B_c$. This is to be achieved by augmenting the laws the controller must satisfy.

Characterize the subsystems of $B_w$ that can be thus achieved.}\\

The controller system $B_c$, and its subsystems obtained by restriction, mediate through the system $B$ (governed by the laws in the submodule $P \subset A^k$) to effect changes in the system $B_w$. This imposes a priori constraints on the possible subsystems of $B_w$ that can be attained by the above process. 

\begin{lemma} The subsystem $\iota_w^{-1}(B)$ of $\pi_w(B)$ is unchanged by additions to the controller laws.
\end{lemma}
\noindent Proof: Imposing additional laws to restrict the behaviour $B_c$ of the controller, translates to specifying an $A$-submodule of $A^{k_c}$ strictly containing the submodule $\iota_c^{-1}(P)$. These laws correspond to laws of the form $(0, q) \in A^{k_w+k_c}$ that are not in $P$. The addition of such laws to $P$ in turn results in restricting the system $B$ to a subsystem. However, the submodule $P'$, generated by $P$ and these new laws, satisfies $\pi_w(P') = \pi_w(P)$. As $\ker_\F(\pi_w(P)) = \iota_w^{-1}(B)$, this subsystem of $B_w$ remains unchanged when $P$ is enlarged to $P'$.  \hspace*{\fill}$\square$

\begin{corollary} The possible subsystems of $\F^{k_w}$ that can be achieved by restricting $\pi_c(B)$ with additional controller laws, all contain $\iota_w^{-1}(B)$, and are contained in $\pi_w(B)$. 
\end{corollary}
\noindent Proof: Together with the above lemma, it suffices to observe that additional controller laws results in a larger set of laws that the system $B$ must satisfy. Let it be given by a submodule $P'$ containing $P$. It follows that $\iota_w^{-1}(P) \subset \iota^{-1}_w(P')$, and hence that the controlled behaviour $\ker_\F(\iota_w^{-1}(P'))$ must be contained in $B_w$. \hspace*{\fill}$\square$\\

Thus, $\iota_w^{-1}(B)$ is a residual subsystem of $B_w$, in the sense that every subsystem of $B_w$ that can be achieved by augmenting the controller $B_c$, contains it. \\

Dual to the above corollary is the following lemma which characterises  the subsystems of $\F^{k_c}$ that can possibly restrict the behaviour of the system $B_w$.

\begin{lemma} Every subsystem of $\F^{k_c}$ that can restrict $\pi_w(B)$ is contained in $\pi_c(B)$, and can be assumed to contain $\iota_c^{-1}(B)$.
\end{lemma}
\noindent Proof: The first containment follows exactly as in the proof of the above corollary.

Next, let $M$ be a submodule of $A^{k_c}$, and let $q \in M \setminus \pi_c(P)$. Then there is no $p \in A^{k_w}$ such that $(p,q) \in P$, and thus $\iota_w^{-1}(P + (0, q)) = \iota_w^{-1}(P)$.  This implies that the addition of the law $q$ to $\iota_c^{-1}(P)$ leaves $B_w$ unchanged. Thus we may assume that $M \subset \pi_c(P)$ by replacing $M$ with $M \cap \pi_c(P)$, and hence that $\iota_c^{-1}(B) \subset \ker_\F(M)$.
\hspace*{\fill}$\square$\\

Corollary 2.1 and Lemma 2.2 above, show that the variables $w$ that are to be controlled, and the control variables $c$, satisfy identical restrictions. In other words, $B_w$ can be restricted to a subsystem containing $\iota_w^{-1}(B)$, by a controller that is contained in $B_c$, and which contains $\iota_c^{-1}(B)$. The above control problem is thus symmetric in the $w$ and $c$ variables. \\

In light of these results, we make the following definition.
\vspace{.25cm}

\noindent Definition: An $A$-submodule of $A^{k_w}$  containing $\iota_w^{-1}(P)$, and which is contained in $\pi_w(P)$, is said to be admissible with respect to $P$ (similarly for submodules  of $A^{k_c}$ containing $\iota_c^{-1}(P)$ and contained in $\pi_c(P)$).

\begin{proposition} Let $\Phi$ assign an $A$-submodule $N$ of $A^{k_w}$ to the submodule $\Phi(N) = \iota_c^{-1}(\iota_w(N) + P)$ of $A^{k_c}$. Then $\Phi$ is a bijection between the admissible submodules of $A^{k_w}$ and the admissible submodules of $A^{k_c}$, with respect to $P$.
\end{proposition}
\noindent Proof:  If $N = \iota_w^{-1}(P)$, then $(\iota_w(N) + P) = P$, and so $\Phi(\iota_w^{-1}(P)) = \iota_c^{-1}(P)$.  

Now let $N = \pi_w(P)$. For every $(p, q) \in P$, $(p, 0)$ is in $\iota_w(N)$, hence $(0, q)$ is in $\iota_w(N) + P$. Thus every $q$ in $\pi_c(P)$ is in $\iota_c^{-1}(\iota_w(N) + P)$, and hence $\Phi(\pi_w(P)) = \pi_c(P)$. 

As the assignment $\Phi$ is inclusion preserving, it follows that it maps an admissible submodule of $A^{k_w}$ to an admissible submodule of $A^{k_c}$, with respect to $P$.

Similarly, for an $A$-submodule $M$ of $A^{k_c}$, define $\Psi (M) = \iota_w^{-1}(\iota_c(M)+ P)$. It assigns admissible submodules of $A^{k_c}$ with respect to $P$ to admissible submodules of $A^{k_w}$. It is easily verified that $\Phi$ and $\Psi$ are inverses of one another, and hence it follows that they are both bijections.    \hspace*{\fill}$\square$\\

We can now characterize the achievable subsystems of $B_w$. For this purpose,
we assume further that $\F$ is a cogenerator. Thus $\F$ could be either $\Dr$ or $\Cinf$ in the case of partial differential operators, or $\C^{\N^n}$ in the case of difference operators.

\begin{theorem}Let $B = \ker_\F(P)$ be the system in $\F^{k_w+k_c}$ defined by the submodule $P \subset A^{k_w+k_c}$. Let the signal space $\F$ be an injective $A$-module, which is also a cogenerator. Then every subsystem of $\pi_w(B)$ containing $\iota_w^{-1}(B)$ can be achieved by a unique controller contained in $\pi_c(B)$ and containing $\iota_c^{-1}(B)$.
\end{theorem}
\noindent Proof: Let $B'$ be a subsystem of $\pi_w(B)$ containing $\iota_w^{-1}(B)$. As $\F$ is an injective cogenerator, $B'$ equals $\ker_\F(N)$, for a unique submodule $N$ of $A^{k_w}$. This submodule $N$ is admissible with respect to $P$, i.e., it satisfies $\iota_w^{-1}(P) \subset N \subset \pi_w(P)$. By the above proposition $\Phi(N)$ equals an admissible submodule of $A^{k_c}$, say $M$.

Suppose that the laws $\iota_c^{-1}(P)$ of the system $\pi_c(B)$ are augmented to this submodule $M$ by the design of a controller. Then the laws of $B$, namely the submodule $P \subset A^k$, are augmented to the submodule $\iota_c(M) + P$. The projection of the resultant system, $\pi_w(\ker_\F(\iota_c(M) + P))$, to $\F^{k_w}$ is a system whose laws are given uniquely by $\i_w^{-1}(\iota_c(M) + P)$, namely $\Psi(M)$ of the above proposition. As $\Psi$ is inverse to $\Phi$, $\Psi(M) = N$.

Thus $B'$ is achieved by the unique subsystem of $\pi_c(B)$ defined by the submodule $M$. \hspace*{\fill}$\square$

\begin{corollary} Amongst all the controllers that restrict $\pi_w(B)$ to a subsystem $B'$, there is a unique minimal one whose laws are derived from the laws $P$ of $B$.
\end{corollary}
\noindent Proof: As $\F$ is an injective cogenerator, the subsystem $B'$ is defined uniquely by a submodule $N \subset A^{k_w}$, admissible with respect to $P$. Consider the submodule $\Phi(N) \subset \pi_c(P)$ determined by the correspondence of Proposition 2.3. By Theorem 2.1, the controller determined by $\Phi(N)$ restricts $B$ to $B'$. The expression for $\Phi$ shows that the laws of this controller are derived from the laws $P$.

Now suppose that $M$ is a set of laws that defines a controller which restricts $B_w$ to $B'$. By Lemma 2.2, it follows that $M' = M \cap \pi_c(P)$ also restricts $B$ to $B'$. Again as $\F$ is an injective cogenerator, it follows that $M'$ must equal $\Phi(N)$.

Thus, the laws that determine any controller which restricts $B_w$ to $B'$, must contain $\Phi(N)$, and therefore $\Phi(N)$ is the unique minimal controller. \hspace*{\fill}$\square$\\

\noindent Remark 2.3: The description of the above minimal controller is an instance of the Internal Model Principle, in the sense that it has sufficient information about the system $B$ built into it, \cite{fw}.\\

Recollect the notion of the canonical controller from \cite{sc}. Given an achievable subsystem $B'$ of $\pi_w(B)$, its canonical controller is the subsystem of $\pi_c(B)$ defined by $C_{can} = \{ f_c \in \F^{k_c} ~| ~ \exists ~f_w \in B' ~ with ~ (f_w,f_c) \in B \}$.

\begin{corollary} The unique minimal controller of the above corollary is the canonical controller of $B'$. 
\end{corollary}
\noindent Proof: Let $B' = \ker_\F(N)$, for a unique submodule $N$ of $A^{k_w}$, admissible with respect to $P$. Then, by definition, the controller $\ker_\F(\Phi(N))$ of the above corollary is obtained from $B$ by restricting $B_w$ to $B'$. This is precisely the canonical controller of van der Schaft. \hspace*{\fill}$\square$\\

These results also answer Problem (ii) above on the structure of controllers.\\

\noindent Remark 2.4: If we set $k = 1$, then we are in the realm of systems defined by ordinary differential operators, and all the above results specialise to the results of \cite{sc} and \cite{wt}.

There are several other issues related to the construction of the canonical controller, especially the notion of regular implementation (for instance \cite{pr1} and \cite{np}), which are, however, not pursued here.

\section{control of the electro-magnetic field}
We study Maxwell's equations in the context of the results of the previous section.\\

Let $A$ be the ring $\C[\del_x, \del_y, \del_z, \del_t]$ of differential operators on space and time. The Maxwell equations (in Gaussian units) are 
\[\nabla \cdot E - 4\pi \rho = 0, \hspace{.3cm}\nabla \cdot B = 0, \]
\[\nabla \times E + \frac{1}{\tt{c}} \del_t B = 0, \hspace{.2cm} \nabla \times B - \frac{1}{\tt{c}}(4\pi J +  \del_t E) = 0 ~, \]
where $E, B$ are the electric and magnetic fields, $\rho, J$, the electric charge and electric current densities, and $\tt{c}$, the speed of light.
 
The partial differential operator defined by  these equations is \[P(\del): \F^{10} \rightarrow \F^8 ~,\] where

{\small 
\[
P(\del)  = \left(
\begin{array}{lccccccccc}
 \phantom{-} \del_x & \del_y & \partial_z & 0 & 0 & 0 & -4 \pi \phantom{-} & 0 & 0 & 0 \\
\phantom{-} 0 & 0 & 0 & \del_x & \del_y & \del_z & 0 & 0 & 0 & 0 \\
\phantom{-} 0 & -\del_z \phantom{-} & \del_y & \frac{1}{\tt{c}}\del_t \phantom{-}& 0 & 0 & 0 & 0 & 0 & 0 \\
\phantom{-} \del_z & 0 & -\del_x \phantom{-} & 0 & \frac{1}{\tt{c}}\del_t \phantom{-}& 0 & 0 & 0 & 0 & 0 \\
-\del_y  & \del_x & 0 & 0 & 0 & \frac{1}{\tt{c}}\del_t \phantom{-}& 0 & 0 & 0 & 0 \\
\frac{1}{\tt{c}} \del_t  & 0 & 0 & 0 & \del_z  & -\del_y \phantom{-} & 0 & \frac{4\pi}{\tt{c}} & 0 & 0 \\
\phantom{-} 0 & \frac{1}{\tt{c}} \del_t \phantom{-}& 0 & -\del_z \phantom{-} & 0 &   \del_x  & 0 & 0 & \frac{4\pi}{\tt{c}} & 0 \\
\phantom{-} 0 & 0 & \frac{1}{\tt{c}}\del_t \phantom{-} &  \del_y & - \del_x \phantom{-} & 0 & 0 & 0 & 0 & \frac{4\pi}{\tt{c}}
\end{array}
\right)
\]}and where $\F$ is either $\Dr(\R^4)$ or $\Cinf(\R^4)$. The eight rows of this matrix correspond to the two equations involving divergence, and the six equations involving curl. It operates on 
$(E_1, E_2, E_3, B_1, B_2, B_3, \rho, J_1, J_2, J_3) \in \F^{10}$, the components of which are the components of $E, B, \rho$ and $J$. The electro-magnetic system is the kernel $\ker_\F(P(\del)) $ of this operator.

The problem is to control the electric and magnetic fields by suitable choices of the control variables $\rho$ and $J$. In the notation of the previous section, we have $w = (E_1, E_2, E_3, B_1, B_2, B_3)$, and $c = (\rho, J_1, J_2, J_3)$.

Let $P \subset A^{10}$, be the submodule generated by the rows of $P(\del)$.  It determines the submodules $\pi_w(P)$ and $\iota_w^{-1}(P)$ of $A^6$, and the submodules $\pi_c(P), ~\iota_c^{-1}(P)$ of $A^4$. An elementary calculation shows that $\iota_c^{-1}(P)$ is the submodule of $A^4$ generated by the continuity equation $\del_t \rho + \nabla \cdot J = 0$, whereas $\pi_c(P) = A^4$. 

Similarly, the submodule $\iota_w^{-1}(P) \subset A^6$ is generated by the laws $\nabla \cdot B = 0$ and $\nabla \times E + \frac{1}{\tt{c}}\del_t B = 0$; these are the `homogeneous' Maxwell equations given by the submatrix defined by rows 2 to 5, and columns 1 to 6 of $P(\del)$. Finally, the submodule $\pi_w(P)$ is generated by the rows of the $8 \times 6$ submatrix of $P(\del)$ defined by its first 6 columns. They are the Maxwell equations in vacuum, namely the homogeneous equations above, together with $\nabla \cdot E = 0$ and $\nabla \times B - \frac{1}{\tt{c}} \del_t E = 0$. 
 
By Lemma 2.2, the laws governing the control variables $\rho$ and $J$ can be any $A$-submodule $M$ of $A^4$ containing the continuity equation. Then, $\rho$ and $J$ would be restricted to lie in $\ker_\F(M)$. In other words, the control variables can be restricted by any system of differential equations containing the continuity equation. Thus, charge and current could be considered to play the classical role of inputs.

By Corollary 2.1, the restriction of $\rho$ and $J$ by the laws in $M$, results in a unique subsystem of the system $\ker_\F(\iota_w^{-1}(P))$ of homogeneous solutions, and containing the vacuum solutions $\ker_\F(\pi_w(P))$. This is the system determined by the admissible submodule $\Psi(M) \subset A^6$ (in the notation of Proposition 2.3). 

In other words, every achievable subsystem of the electro-magnetic field lies between two systems, one, the solutions of the homogeneous Maxwell equations, and the other, the solutions of the vacuum equations. They are
obtained by imposing additional differential constraints on the current and charge densities, and these constraints translate to laws that the electric and magnetic fields must satisfy, in addition to the homogeneous equations. The canonical controller is determined by  the single criterion that the solutions of the controller equations contain the solutions of the continuity equations. As every controller must satisfy the continuity equation, it folows that there is only one controller that accomplishes a given restriction, and hence that this controller is the canonical one.\\

As an example, suppose that the electric charge density $\rho$ is set to 0, by the imposition of the law defined by the cyclic submodule of $A^4$ generated by $(1, 0, 0, 0)$. Let $M$ be the submodule of $A^4$ generated by this law together with the continuity equation. Thus, suppose that $\rho = 0$, and hence that $\nabla \cdot J = 0$. Then $\Phi(M) \subset A^6$ is the submodule  generated by the homogeneous equations together with $\nabla \cdot E = 0$, and the system $B_w$ is restricted to $\ker_\F(\Phi(M))$ by this control action. \\

Conversely, by Theorem 2.1, every electro-magnetic system contained between these two extremes is achievable by suitably restricting  electric charge and electric current, in addition to satisfying the continuity equation. 
 
This is precisely the physics of the electro-magnetic field. \\

\section{Acknowledgement} I am grateful to Madhu Belur for explaining to me the notion of the canonical controller. I  thank Virendra Sule for several discussions on the Internal Model Principle, and Alok Laddha for help with Gaussian units.\\

\end{document}